\documentclass[11pt,reqno]{amsart}
\usepackage{euscript,amssymb}
\usepackage[arrow,matrix,curve]{xy}
\usepackage{array,longtable}
\usepackage{a4wide}

\usepackage{enumitem}

\newenvironment{m-theorem}{%
\vskip4pt\refstepcounter{stff}\trivlist \itemindent 0pt
\item[\hskip\labelsep\bf Theorem \thestff]%
\it\ignorespaces}{\endtrivlist\vskip4pt}%
\newenvironment{m-proposition}{%
\vskip4pt\refstepcounter{stff}\trivlist \itemindent 0pt
\item[\hskip\labelsep\bf Proposition \thestff]%
\it\ignorespaces}{\endtrivlist\vskip4pt}%
\newenvironment{m-corollary}{%
\vskip4pt\refstepcounter{stff}\trivlist \itemindent 0pt
\item[\hskip\labelsep\bf Corollary \thestff]%
\it\ignorespaces}{\endtrivlist\vskip4pt}%
\newenvironment{m-lemma}{%
\vskip4pt\refstepcounter{stff}\trivlist \itemindent 0pt
\item[\hskip\labelsep\bf Lemma \thestff]%
\it\ignorespaces}{\endtrivlist\vskip4pt}%
\newenvironment{m-definition}{%
\vskip4pt\refstepcounter{stff}\trivlist \itemindent 0pt
\item[\hskip\labelsep\bf Definition \thestff]%
\ignorespaces}{\endtrivlist\vskip4pt}%
\newenvironment{m-notation}{%
\vskip4pt\refstepcounter{stff}\trivlist \itemindent 0pt
\item[\hskip\labelsep\bf Notation \thestff]%
\ignorespaces}{\endtrivlist\vskip4pt}%
\newenvironment{m-example}{%
\vskip4pt\refstepcounter{stff}\trivlist \itemindent 0pt
\item[\hskip\labelsep\bf Example \thestff]%
\ignorespaces}{\endtrivlist\vskip4pt}
\newenvironment{m-remark}{%
\vskip4pt\refstepcounter{stff}\trivlist \itemindent 0pt
\item[\hskip\labelsep\bf Remark \thestff]%
\ignorespaces}{\endtrivlist\vskip4pt}
\newenvironment{m-question}{%
\vskip4pt\refstepcounter{stff}\trivlist \itemindent 0pt
\item[\hskip\labelsep\bf Question.]%
\ignorespaces}{\endtrivlist\vskip4pt}%
\newenvironment{thm-nono}[1]{
\vskip5pt\trivlist \itemindent 0pt
\item[\hskip\labelsep\bf Theorem #1.]%
\it\ignorespaces}{\endtrivlist\vskip5pt}%
\newenvironment{m-thank}{%
\vskip4pt\trivlist \itemindent 0pt
\item[\hskip\labelsep\it Acknowledgments]%
\ignorespaces}{\endtrivlist\vskip4pt}%
\newenvironment{m-proof}{%
\vskip4pt\trivlist \itemindent 0pt
\item[\hskip\labelsep\it Proof.]%
\ignorespaces}{\hfill$\Box$\endtrivlist\vskip4pt}%
\let\euf\EuScript 
\let\cal\mathcal
\let\mbb\mathbb

\let\mfrak\mathfrak
 
\DeclareFontFamily{OT1}{rsfs}{}
\DeclareFontShape{OT1}{rsfs}{n}{it}{<->rsfs10}{}
\DeclareMathAlphabet{\crl}{OT1}{rsfs}{n}{it}
\newcommand\uset[2]{{\disp\mathop{\mbox{$#2$}}_{#1}}}

\newcommand\ouset[3]{{\overset{#2}{\underset{#1}#3}}}
\let\ovl\overline
\let\unl\underline

\numberwithin{equation}{section}
\numberwithin{figure}{section} 
\newcommand\Aut{\operatorname{\textrm{Aut}\kern1pt}}
\newcommand\End{\operatorname{\rm{End}\kern1pt}}
\newcommand\cEnd{\operatorname{\mathcal{E}\kern-1pt\textit{nd}\kern1pt}}
\newcommand\Ext{\mathop{\rm Ext}\nolimits}

\newcommand\cHom{\operatorname{\mathcal{H}\kern-1pt\textit{om}\kern1pt}}

\newcommand\Ker{{\rm Ker}}

\newcommand\Pic{\mathop{\rm Pic}\nolimits}

\newcommand\eA{{\euf A}}
\newcommand\eF{{\euf F}}
\newcommand\eG{{\euf G}}
\newcommand\eI{{\euf I}}
\newcommand\eO{{\euf O}}
\newcommand\eL{{\euf L}}
\newcommand\cU{{\cal U}}
\newcommand\ball{{\mbb B}}
\newcommand\bone{{1\kern-0.57ex\rm l}}
\newcommand{\bloc}{\mathop{\rm base\,locus}\nolimits}
\newcommand{\sbloc}{\mathop{\rm stable\,base\,locus}\nolimits}
\newcommand\cd{{\rm cd}}
\newcommand\codim{{\rm codim}}
\let\dashto\dashrightarrow
\let\del\partial
\let\disp\displaystyle
\newcommand{\dvs}{\mathop{\rm divisor}\nolimits}
\let\ges\geqslant
\newcommand\Gl{{\rm GL}}
\let\hra\hookrightarrow

\newcommand{\kk}{{\mbb C}}
\let\les\leqslant
\newcommand\lran[1]{{\langle #1\rangle}}
\newcommand\mx{{\rm max}}

\let\nit\noindent
\newcommand\res{\mathop{\rm res}\nolimits}
\newcommand\rd{{\rm d}} 
\newcommand\rk{{\rm rk}}
\let\si\sigma
\let\Si\Sigma
\let\sm\setminus
\newcommand\Spec{\mathop{\rm Spec}\nolimits}
\let\srel\stackrel
\let\surj\twoheadrightarrow
\let\tld\tilde
\let\unl\underline
\let\veps\varepsilon
\let\vphi\varphi

\author{Mihai Halic}
\email{mihai.halic@gmail.com}


\begin{document}

\title[Splittings induced by restrictions to divisors]%
{Splitting criteria for vector bundles induced by restrictions to divisors}

\begin{abstract}
In this article we deduce criteria for the splitting and the triviality of vector bundles, by restricting them to partially ample divisors. This allows to study the problem of splitting on the total space of fibre bundles. The statements are illustrated with a number of examples. For products of minuscule homogeneous varieties, our results allow to test the splitting of vector bundles by restricting them to products of Schubert $2$-planes. 

The triviality criteria obtained inhere are particularly suited to Frobenius split varieties, whose splitting is defined by a section in the anti-canonical line bundle. As an application, we prove that a vector bundle on a smooth toric variety $X$, whose anti-canonical bundle has stable base locus of co-dimension at least three, is trivial when its restrictions to the invariant divisors are trivial, with trivializations compatible along the various intersections.
\end{abstract}

\maketitle


\section*{Introduction}

Although the problem of deciding the splitting of vector bundles is very classical, only relatively few cases have been settled: there are cohomological criteria for products of projective spaces and quadrics \cite{c-mr,ba-ma}, Grassmannians \cite{ottv,mals}, hypersurfaces in projective spaces \cite{swd,bis+rav}. Splitting criteria corresponding to restrictions are useful because they yield dimensional reductions: the problem is reduced to a (usually much) lower dimensional variety, where one can use further cohomological tools. 

Horrocks' criterion \cite{hor} states that a vector bundle on $\mbb P^n$, $n\ges 3$, splits if and only if its restriction to some hyperplane $D\cong\mbb P^{n-1}$ splits. This was generalized in \cite{ba}, where the author restricted vector bundles on `\emph{Horrocks varieties}' to \emph{ample divisors}. The ampleness assumption excludes several natural situations, \textit{e.g.} the case of  morphisms, where one wishes to restrict vector bundles either to pre-images of ample divisors or to relatively ample ones.  Also, the `Horrocks variety' assumption is a rather restrictive cohomological property.

The goal of this note is to generalize the splitting criterion in \textit{op. cit.} to include $q$-ample divisors; this covers the case of morphisms mentioned before. We obtain two types of results: \emph{splitting} and \emph{triviality} criteria. 

\begin{thm-nono}{(splitting criteria)}  
Let $(X,\eO_X(1))$ be a smooth, complex projective variety, with $\dim X\ges3$. Let $\crl V$ be a vector bundle on $X$, $\crl E:=\cEnd(\crl V)$ the bundle of endomorphisms, $\eL\in\Pic(X)$ be $q$-ample, and $D\in|d\eL|$. The equivalence 
$$\;[\;\crl V\;\text{splits}\;\Leftrightarrow\;\crl V\otimes\eO_D\;\text{splits}\;]\;$$ 
holds in any of the following cases: 
\begin{enumerate}
\item[(a)] 
If $q\les\dim X-3$, $D$ is either reduced or irreducible, and $H^1(\crl E_D\otimes\eL_D^{-a})=0$ for all $a\ges d$.  The parameter $d$ is bounded from below by a linear function in the regularity of $\crl E$ with respect to $\eO_X(1)$. 
\smallskip
\item[(b)] 
If $X$ is $2$-split, $q\les\dim X-4$, and $D$ is smooth.
\smallskip
\item[(c)] 
If $X$ is $1$-split, $q\les\dim X-4$, $\eL$ is globally generated, and $D\in|\eL|$ is very general.
\end{enumerate} 
\end{thm-nono}
The conditions $1$-, $2$-split are respectively the notions of `splitting' and `Horrocks variety' in \cite{ba}. Bakhtary's result corresponds to the case (b) above with $q=0$. 
The statement (a) is more effective in the case where $\eL$ is relatively ample with respect to a morphism (cf. Theorem \ref{thm:rel-ample}): it suffices $D\in|\eL|$ to be weakly normal and $\crl E\otimes\eL$ be relatively ample. 
\smallskip

We illustrate the advantage of allowing $q$-ample line bundles by discussing several explicit examples. First we simplify the cohomological splitting criteria for products of projective spaces and quadrics \cite{c-mr,ba-ma}. They involve a large number of cohomological tests and, by restricting to `sub-products', the number of the tests is massively reduced. 

Second, we use the results in \cite{hal2} and we deduce a splitting criterion for vector bundles on varieties $X$ which are products of minuscule homogeneous varieties (cf. Theorem \ref{thm:minuscule}). To our best knowledge, currently there are \emph{no known results} in this direction; even a product of Grassmannians seems to be uncovered. We prove that, for such $X$, it is enough to test the splitting on appropriate products of planes $\mbb P^2$ (Schubert subvarieties) in $X$.  \smallskip

The trivializable vector bundles are particular cases of the split ones, so the triviality criteria below hold in greater generality. Notably, one can eliminate the conditions $1$-, $2$-split. 

\begin{thm-nono}{(triviality criteria)} 
Let $X,\eL,\crl V$ be as above, and $D\in|\eL|$. 

\nit The equivalence 
$\;[\,\crl V\;\text{is trivial}\;\Leftrightarrow\;\crl V_D\;\text{is trivial}\,]\;$ 
holds in any of the following cases: 
\begin{enumerate}
\item[(a)]
If $\eL\in\Pic(X)$ is semi-ample and $(\dim X-3)$-ample.  
\item[(b)] 
If $\eL$ is relatively ample for a morphism $X\srel{f}{\to} Y$ of relative dimension at least three. 
\item[(c)] 
If the anti-canonical bundle $\omega_X^{-1}$ is $(\dim X-3)$-ample, $X$ is Frobenius split by a power of a section $\si$ in $\omega_X^{-1}$, and $D=\dvs(\si)$. 
\end{enumerate}
\end{thm-nono}
The condition (c) above is particularly suited for spherical varieties (\textit{e.g.} toric varieties), because they satisfy the assumption about the Frobenius splitting. We elaborate on the case of toric varieties. 

\begin{thm-nono}{\kern-.75ex}
Let $X$ be a smooth toric variety and $\Delta$ be its boundary divisor. We assume that 
$$
\codim\big(\sbloc(\omega_{X}^{-1})\,\big) \ges 3.
$$ 
Then, for a vector bundle $\crl V$ on $X$, one has the equivalences: 
\begin{enumerate}
\item[(a)] 
$[\,\crl V\mbox{ splits}\;\Leftrightarrow\;\crl V_{\Delta_m}\mbox{ splits}\,]$, for $m\gg0$. 
\item[(b)] 
$[\,\crl V\mbox{ is trivial}\;\Leftrightarrow\;\crl V_{\Delta}\mbox{ is trivial}\,]$.
\end{enumerate}
\end{thm-nono}

The splitting criteria obtained in this article are based on two technical ingredients: the `universal' criterion \ref{prop:h1}, on one hand, and various Kodaira-type vanishing theorems for $q$-ample line bundles, on the other hand. 

For this reason, the role of the appendices is twofold: first, to recall the definitions and properties of the $q$-ampleness (cf. \cite{totaro,ottem}) and Frobenius splitting (cf. \cite{brion-kumar}) which are used in the body of the article; second, they contain a few (possibly) noteworthy results: 
\begin{itemize}[leftmargin=*]
\item 
a Kodaira vanishing theorem for relatively ample line bundles on weakly normal varieties (cf. \ref{thm:global-gen}(iii)) and for $q$-ample line bundles on Frobenius-split varieties (cf. \ref{thm:F-split});
\item 
a Picard-Lefschetz property in the relative setting (cf. the Theorem \ref{thm:pic-iso}); 
\item 
a $q$-ampleness criterion for line bundles which are not necessarily globally generated (cf. Theorem \ref{thm:q-ample}). 
\end{itemize}
Throughout this article, $X$ stands for a smooth projective variety over $\kk$ of dimension at least three.


\section{The general splitting principle}\label{sct:split-vb}

\begin{m-definition}\label{def:split}
Let $T$ be a scheme defined over $\kk$ and $S$ be a closed subscheme of it; we assume that 
$H^0(\eO_S)=H^0(\eO_T)=\kk.$ 
For a locally free sheaf (a vector bundle) $\crl V_T$ of rank $r$ on $T$, we denote $\crl E_T:=\cEnd(\crl V_T)$ the sheaf of endomorphisms; let $\crl V_S:=\crl V_T\otimes_{\eO_T}\eO_S$, \textit{etc}.  

An \emph{eigenvalue} of $h_S\in H^0(\crl E_S)$ is a complex root of the polynomial 
\begin{equation}\label{eq:ph}
p_{h_S}:=\det\big(t\bone-h_S\big)\in
 H^0\big(\,
\cEnd(\det\crl V_S)
\,\big)[t]
= H^0(\eO_S)[t]=\kk[t].
\end{equation}
We say that $\crl V_T$ \textit{splits} if it is isomorphic to a direct sum of $r$ invertible sheaves (line bundles) on $T$. 
\end{m-definition}

Let $S,h_S$ be as above. If $\veps\in\kk$ is an eigenvalue of $p_{h_S}$, then 
$\Ker(\veps\bone-h_S)\subset\crl V_S$ is a non-zero $\eO_S$-module: indeed, for a closed point $x\in S_{\rm red}\subset S$ with maximal ideal $\mfrak m_x\subset\eO_S$, $\veps$ is an (usual) eigenvalue of $h_S\otimes\frac{\eO_S}{\mfrak m_x}\in\End\big( \crl V_S\otimes\frac{\eO_S}{\mfrak m_x} \big)$. 

\begin{m-lemma}\label{lm:split}
Let the notation be as above. Then the following statements hold:
\begin{enumerate}
\item 
$\crl V_S$ splits if and only if there is $h_S\in H^0(\crl E_S)$ with $r$ pairwise distinct eigenvalues. 
\item  
If $H^0(\crl E_{T})\to H^0(\crl E_{S})$ is surjective, then $\crl V_{T}$ splits if and only if $\crl V_S$ splits.
\end{enumerate}
\end{m-lemma}

\begin {m-proof}
(i) 
Suppose that $h_S\in H^0(\crl E_S)$ has pairwise distinct eigenvalues $\veps_1,\ldots,\veps_r\in\kk$, and let $\ell_j:=\Ker(\veps_j\bone-h_S)\neq0$. The (polynomial) identity 
\\ \centerline{$
\mbox{$
1=\ouset{j=1}{r}{\sum}c_jp_j(t),\quad\text{with}\;
c_j:=\biggl(\underset{k\neq j}{\prod}(\veps_j-\veps_k)\biggr)^{-1}\in\kk,
\quad p_j(t):=\underset{k\neq j}{\prod}(t-\veps_k)\in\kk[t],
$}
$}\\[.5ex] 
implies $\crl V_S\,{=}\ouset{j=1}{r}{\sum}\,{\rm Image}\big(p_j(h_S)\big)$. 
Moreover, the Cayley-Hamilton theorem yields: 
\\ \centerline{$
{\rm Image}\big(p_j(h_S)\big)\subset\ell_j\;\Rightarrow\; 
\mbox{$\crl V_S\,{=}\ouset{j=1}{r}{\sum}\,\ell_j.$}
$}\\[1ex]
We claim that the sum is direct. Indeed, if $v_j\in\ell_j$ for all $j$, then holds: 
$$
\begin{array}{rcll}
v_1+v_2+\ldots+v_r=0 &{\Rightarrow}& p_1(h_S)(v_1)=0, 
&
\mbox{(since $t-\veps_j\mid p_1(t)$, for $j\ges 2$),}
\\[1ex]
v_1=\ouset{j=2}{r}{\sum}c_j\cdot p_j(h_S)(v_1)&=&0,\;etc.
&
\mbox{(since $v_1\in\Ker(\veps_1\bone-h_S)$).}
\end{array}
$$
Hence $\ell_j$, $j=1,\ldots,r$, are (locally) projective $\eO_S$-modules, so they are locally free (cf. \cite[Ch.\,II, \S 5.2, Th\'eor\`eme 1]{bour-comm1-4}) of rank one. 

\nit(ii) 
If $\crl V_S$ splits, there is $h_S\in H^0(\crl E_S)$ with $r$ pairwise distinct complex eigenvalues; take an extension $h_{T}\in H^0(\crl E_{T})$ of it. The equation \eqref{eq:ph} shows that $p_{h_T}{=}\,\det(t\bone-h_{T})\in\kk[t]$, so $p_{h_T}=p_{h_S}$, hence $h_T$ has the same eigenvalues as $h_S$. 
\end {m-proof}

\begin{m-definition}\label{def:D-thick}
Let $\eL$ be an invertible sheaf (a line bundle) on $X$ and $D\in|d\eL|$ be an effective divisor. For $m\ges 0$, the $m$-th order thickening $D_m$ of $D$ is the subscheme of $X$ defined by the ideal $\eI_D^{m+1}$, where $\eI_D=\eO_X(-D)\cong\eL^{-d}$. 
\end{m-definition}
The structure sheaves of successive thickenings fit into the exact sequences: 
\begin{equation}\label{eqn:thick}
0\to\eL_D^{-dm}\to\eO_{D_m}\to\eO_{D_{m-1}}\to0,\quad m\ges1.
\end{equation}
We will apply the Lemma \ref{lm:split} mostly in the case $T=X, S=D_m$, for suitable $m$.

\begin{m-remark}
If $D$ is an effective divisor of $X$, the surjectivity of $H^0(\crl E_X)\to H^0(\crl E_D)$ is implied by the vanishing of $H^1(\eO_X(-D)\otimes\crl E_X)$. Thus, at a certain extent, the $(\dim X-2)$-amplitude of $\eO_X(D)$ is the weakest possible assumption which allows to deduce splitting criteria for vector bundles by restricting them to $D$. 
\end{m-remark}

In the framework of formal schemes, we have the following very general statement.

\begin{m-proposition}\label{prop:formal}
Let $D\subset X$ be an effective divisor, $\dim X\ges2$, and $\hat X:=\disp\varprojlim_{m}D_m$ denote the formal completion of $X$ along $D$. If the cohomological dimension $\cd(X\sm D)\les\dim X{-}2$, then $\crl V$ splits if and only if $\crl V\otimes\eO_{\hat X}$ does. The assumption is satisfied if $D$ is $(\dim X{-}2)$-ample.
\end{m-proposition}

\begin{m-proof}
By \cite[Theorem 3.4]{hart-as}, the restrictions 
$$
\kk\cong H^0(\eO_X)\to H^0(\eO_{\hat X}),\quad 
H^0(\crl E)\to H^0(\crl E_{\hat X})
$$
are both isomorphisms; we conclude by \ref{lm:split}(ii). The last claim is \cite[Proposition 5.1]{ottem}.
\end{m-proof}

\begin{m-lemma}\label{lm:q-split}
Let $\eL$ be a $q$-ample line bundle on $X$, with $q\les\dim X-2$. Consider $D\in|\eL|$ which is either reduced or irreducible. Then the following statements hold: 
\begin{enumerate}
\item 
$H^0(\eO_{D_m})=\kk$, for all $m\ges0$; 
\item 
$\crl V_X$ splits if and only if its restriction $\crl V_{D_m}$ splits, for some $m\gg 0$.  
\end{enumerate}
The conclusions (i), (ii) above hold for an arbitrary $D\in|\eL|$ in any of the situations enumerated below:
\begin{equation}\label{eq:C}
\mbox{
\begin{minipage}{0.85\linewidth}
\nit{\rm(a)} 
$\eL$ is semi-ample;

\nit{\rm(b)} 
$\eL$ is relatively ample for a morphism $X\srel{f}{\to}Y$, $\dim X-\dim Y\ges2$;

\nit{\rm(c)} 
$X$ is an F-split variety (cf. appendix \ref{sct:frob}).
\end{minipage}
}
\end{equation}
\end{m-lemma}

\begin {m-proof}
(i) The $(\dim X-2)$-amplitude of $\eL$ implies that $\kk= H^0(\eO_X)\to H^0(\eO_{D_m})$ is an isomorphism for $m\gg0$, so $D$ is connected. 

If $D$ is reduced, then we have $H^0(\eO_D)=\kk$. As $H^0(\eL_D^{-m})=0$ for $m\gg0$, the same holds for all $m\ges1$.  The conclusion follows from \eqref{eqn:thick}. 
If $D$ is irreducible, then $D=m_0D_{\rm red}$ for some $m_0\ges1$, and the previous argument shows that $H^0(\eO_{D_{{\rm red},m}})=\kk$ for all $m\ges1$. 

Concerning the final claim, in the cases \eqref{eq:C} we have $H^j(X,\eL^{-m})=0$, for $j=0,1$ and all $m\ges 1$ (cf. \ref{thm:global-gen} and \ref{thm:F-split}), which implies that $H^0(\eO_{D_m})=\kk$ for all $m\ges 0$.

\nit(ii) We apply the Lemma \ref{lm:split}: $H^0(\crl E_X)\to H^0(\crl E_{D_m})$ is surjective if $H^1(\crl E_X\otimes\eL^{-d(m+1)})=0$. This is indeed the case, for large $m$. 
\end {m-proof}

The following general splitting principle, corresponding to restrictions to partially ample divisors, is the root of the results obtained in this article. 

\begin{m-proposition}\label{prop:h1}
Let $\eL\in\Pic(X)$ be $q$-ample, with $q\les\dim X-2$, and $D\in|d\eL|$ be an effective divisor which is either \emph{reduced} or \emph{irreducible}. Assume that 
\begin{equation}\label{eq:h1}
H^1(D,\crl E_D\otimes\eL_D^{-a})=0,\quad\forall a\ges c.
\end{equation}
Then the following properties hold: 
\begin{enumerate}
\item 
$H^1\bigl(X,\crl E\otimes\eL^{-a}\bigr)=0$, for all $a\ges c$. 
\item 
If $d\ges c$ and $\crl V_D$ splits, then $\crl V$ splits too. 
\end{enumerate}
In any of the cases \eqref{eq:C}, $D$ can be \emph{arbitrary}.
\end{m-proposition}

\begin {m-proof}
(i) Denote 
$$a_0:=\max\{a\mid H^1(X,\crl E\otimes\eL^{-a})\neq0\}<\infty.$$ 
The exact sequence $0\to\eL^{-d}\to\eO_X\to\eO_D\to0$ yields 
$$
\ldots\to 
H^1\bigl(\crl E\otimes\eL^{-d-a_0}\bigr)
\to 
H^1\bigl(\crl E\otimes\eL^{-a_0}\bigr)\to
H^1\bigl(\crl E_D\otimes\eL_D^{-a_0}\bigr)
\to\ldots,
$$
with $-d-a_0\les-(a_0+1)$, so the leftmost term vanishes. If $a_0\ges c$, then the rightmost and the middle terms vanish too. This contradicts the definition of $a_0$, hence $a_0<c$.  

\nit (ii) We have $H^0(\eO_D)=\kk$, by \ref{lm:q-split}; since $d\ges c$, $H^0(\crl E)\to H^0(\crl E_D)$ is surjective. 
\end {m-proof}

\begin{m-remark}\label{rmk:c}
(i) The uniform $q$-ampleness property \cite[Theorem 6.4]{totaro} implies that there is a linear function $l(r)=\lambda r+\mu$, with $\lambda,\mu$ depending only on $\eL$, such that \ref{prop:h1}(i) holds for all $a\ges l\big({\rm reg}(\crl E)\big)$, where ${\rm reg}(\crl E)$ stands for the regularity of $\crl E$ with respect to a (fixed) ample line bundle $\eO_X(1)$. 

\nit(ii) The condition \eqref{eq:h1} involves \emph{only} the restriction of $\crl E$ to $D$, which splits by assumption. This feature is helpful because it is easier to decide the vanishing of the cohomology of line bundles, rather than of vector bundles (\textit{i.e.} \ref{prop:h1}(i)). Also, for $q\les\dim X-3$, the condition \eqref{eq:h1} is indeed fulfilled for $c\gg0$. 

\nit(iii) There are two important classes of $q$-ample line bundles: the relatively ample and the pulls-back of ample line bundles with respect to a morphism. We will constantly elaborate on these two situations; the case of a pull-back typically requires stronger hypotheses.
\end{m-remark}

\begin{m-theorem}\label{thm:rel-ample}
Let $\crl V$ be a vector bundle on $X$, $\crl E:=\cEnd(\crl V)$, and 
$f:X\to Y$ be a surjective morphism with $Y$ projective, 
$$
\dim X-\dim Y\ges 3,\;\eL\in\Pic(X) \text{ $f$-relatively ample}.
$$ 
Let $D\in|\eL|$ be a \emph{reduced}, \emph{weakly normal} divisor, and assume that $\crl E_D\otimes\eL_D$ is relatively ample with respect to $D\to Y$. (In particular, it suffices $\crl E\otimes\eL$ to be relatively ample.) Then we have the equivalence:  
$[\,\crl V\;\mbox{splits} \;\Leftrightarrow\; \crl V_D\;\mbox{splits}\,].$
\end{m-theorem}
The weak normality condition for a divisor (cf. \cite[Proposition 4.1]{le-vi}) is explicit, but it is somewhat technical. However, one can see that the condition is satisfied in the following fairly general situation (a particular case of the $WN1$-property \cite[Definition 3.2]{cgm1}): 
\begin{itemize}[leftmargin=*]
\item 
$D=D_1+\dots+D_t$ is reduced, and $D$ is normal away from its self-intersections of a single irreducible component or of two different components; 
\item 
For any point $p\in D$, the local equations (in the analytic topology) of the components of $D$, which are passing through $p$, form a regular sequence in $\hat\eO_{X,p}\cong\kk\{\xi_1,\dots,\xi_{\dim X}\}$. Hence the locus of the points which belong to at least three branches of $D$ have codimension at least two in $D$. 
\item 
$D$ has \emph{generically} normal crossings: at the generic intersection point between two local (analytic) branches of $D$ there are local (analytic) coordinates $\{\xi_1,\dots,\xi_{\dim X}\}$ in $\hat{\eO}_{X,p}$ such that the germ of $D$ at $p$ is given by $\{\xi_1\xi_2=0\}$. 
\end{itemize}

\begin {m-proof} 
Let $\ell\in\Pic(D)$ be a direct summand of $\crl E_D$; by hypothesis, $\ell_a:=\ell\otimes\eL^a$ is relatively ample, for all $a\ges 1$. Note that $D$ is Gorenstein, so it satisfies Serre's condition $S_2$. Then the Theorem \ref{thm:global-gen}(iii) implies that $H^1(D,\ell_a^{-1})=0$, that is the condition \eqref{eq:h1} is fulfilled and we may apply the Proposition \ref{prop:h1}. 
\end {m-proof}


\section{Splitting along divisors: a `deterministic' approach} \label{sct:q-ample}

\begin{m-definition}\label{def:s-split}
For $s\ges 1$, we say that a scheme $S$ is $s$-split if 
\begin{equation}\label{eq:s-split}\tag{$s$-split}
H^j(S,\ell)=0,\text{ for }j=1,\ldots,s,\;\forall\,\ell\in\Pic(S).
\end{equation}
\end{m-definition}

\begin{m-remark}\label{rmk:s-split}
For $s=1,2$ one gets respectively the `splitting' and `Horrocks scheme' notions introduced in \cite{ba}.  Examples of projective varieties satisfying \eqref{eq:s-split} are as follows:
\smallskip
 
\nit(i) arithmetically Cohen-Macaulay varieties $X$---\textit{e.g.} homogeneous spaces, complete intersections in them---with cyclic Picard group (where $s=\dim X-1$), and products of them;

\nit(ii) projective bundles: if $Y$ is $s$-split, and $\euf M_1,\ldots,\euf M_r$, $r\ges s+2$, are line bundles on $Y$, then $X:=\mbb P(\euf M_1\oplus\ldots\oplus\euf M_r)$ is also $s$-split (cf. \cite[Example 4.9]{ba}).
\end{m-remark}

\begin{m-proposition}\label{prop:s-1s}
Let $D$ be an effective $q$-ample divisor on $X$. The following statements hold:

\begin{enumerate}
\item 
If $D$ is $s$-split, with $s\les\dim X-(q+1)$, then $X$ is $s$-split. 
\item 
If $D$ is $s$-split, with $q\les s\les\dim X-(q+1)$, then $X$ is $(s+1)$-split. 
\item 
If $X$ is $(s+1)$-split, then  $D$ is $s$-split in the following cases: 
\begin{enumerate}
\item[(a)] $D$ is smooth and $\eO_X(D)$ is $(\dim X-4)$-positive;
\item[(b)] $D$ is arbitrary, relatively ample for a morphism $X\to Y$, with $Y$ projective and $\dim X-\dim Y\ges4$.
\end{enumerate}
\end{enumerate}
\end{m-proposition}

\begin {m-proof}
(i) We consider the exact sequences 
$$
0\to\eO_X((k-1)D)\to\eO_X(kD)\to\eO_D(kD)\to 0,\;\forall\,k\in\mbb Z,
$$
and tensor them by $\eL\in\Pic(X)$. We obtain surjective homomorphisms 
$H^i(X,\eL((k-1)D))\surj H^i(X,\eL(kD))$, for $i\les s$. 
The $q$-ampleness of $D$ implies that $H^i(X,\eL(kD))=0$ for $k\ll0$, which yields $H^i(X,\eL)=0$. 

\nit(ii) We should verify only that $H^{s+1}(X,\eL)=0$. The previous exact sequence yields 
inclusions $H^{s+1}(X,\eL((k-1)D))\subset H^{s+1}(X,\eL(kD))$, for all $k\in\mbb Z$. Again, the $q$-ampleness of $D$ implies $H^{s+1}(X,\eL(kD))=0$, for $k\gg0$, so $H^{s+1}(X,\eL)=0$. 

\nit(iii) In both cases $\Pic(X)\to\Pic(D)$ is an isomorphism, by the Theorem \ref{thm:pic-iso}. Thus, for any $\ell\in\Pic(D)$ there is $\tld\ell\in\Pic(X)$ such that $\tld\ell_D=\ell$. It remains to take the cohomology of the sequence $0\to\tld\ell(-D)\to\tld\ell\to\ell\to0$. 
\end {m-proof}

\begin{m-theorem}\label{thm:1} 
Let $D$ be an effective divisor on $X$ which is either \emph{reduced} or \emph{irreducible}; 
in the cases enumerated below, $D$ can be \emph{arbitrary}: 
\begin{enumerate}
\item[(a)] 
$\eO_X(D)$ is semi-ample;
\item[(b)]
$D$ is relatively ample for a morphism $X\to Y$, with $Y$ projective, $\dim X-\dim Y\ges2$;
\item[(c)] 
$X$ is an F-split variety (cf. appendix \ref{sct:frob}).
\end{enumerate}
Assume, moreover, that $D$ is $1$-split and $\eO_X(D)$ is $(\dim X-2)$-ample. Then we have the equivalence: 
$
\,[\crl V\;\text{splits}\;\Leftrightarrow\;\crl V_D\;\text{splits}\,].
$ 
\end{m-theorem}

\begin {m-proof}
Since $\crl V_D$ splits and $D$ is $1$-split, it holds $H^1(\crl E_D\otimes\eL_D^{-a})=0$ for all $a\ges 1$. The conclusion follows from the Proposition \ref{prop:h1}. 
\end {m-proof}
The interest in allowing partial ampleness for line bundles, which is considerably weaker than amplitude, is to apply the result for morphisms (\textit{e.g.} fibre bundles).

\begin{m-corollary}\label{cor:relative}
Let $X$ be a smooth, $2$-split, projective variety and $D$ be an effective divisor on it; let $X\srel{f}{\to} Y$ be a morphism. Then the splitting of $\crl V_D$ implies the splitting of $\crl V$ in any of the following cases: 
\begin{enumerate}
\item[(a)] 
$f$ is smooth, $D=f^{-1}(D_Y)$, with $D_Y\subset Y$ a smooth $(\dim Y-4)$-positive divisor;
\item[(b)] 
$D$ is arbitrary, $f$-relatively ample, and $\dim X{-}\dim Y\ges4$.
\end{enumerate}
\end{m-corollary}

\begin {m-proof}
In both situations, the Proposition \ref{prop:s-1s}(iii) implies that $D$ is $1$-split. 
\end {m-proof}

Note that (b) above generalizes \cite[Corollary 4.14]{ba} to the relative case. The result seems to be new even in the case when $X$ is a product (cf. \ref{cor:pxq} below). 

\begin{m-example}\textbf{(Projective bundles).}\label{expl:proj-bdl}
Let $Y$ be a smooth, projective, $1$-split variety. Consider $\euf M,\euf M_1,\ldots,\euf M_t\in\Pic(Y)$, $t\ges3$, and define 
$$X:=\mbb P(\eO_Y\oplus\euf M_1\oplus\ldots\oplus\euf M_t)\srel{f}{\to}Y,\;\; 
D:=\mbb P(\euf M_1\oplus\ldots\oplus\euf M_t)\in|\eO_f(1)|.$$  
($D$ is $1$-split, by \ref{rmk:s-split}(ii).) For any vector bundle $\crl V$ on $X$, we have: [\,$\crl V$ splits $\Leftrightarrow$ $\crl V_D$ splits\,].

By repeatedly applying this method, one reduces the verification of the splitting of $\crl V$ to a $\mbb P^2$-sub-bundle of $X$ over $Y$.
\end{m-example}

\nit\textbf{(Vector bundles on products of projective spaces and quadrics).} 
A splitting criterion for vector bundles on $X_1:=\mbb P^{n_1}\times\ldots\times\mbb P^{n_t}$ is obtained in \cite[Theorem 4.7]{c-mr}. It generalizes Horrocks' criterion, and involves the vanishing of $(n_1+1)\cdot\ldots\cdot(n_t+1)$ cohomology groups. 

The result has been extended in \cite[Theorem 2.14, 2.15]{ba-ma} to products $X_1\times X_2$, where $X_1$ is as above and $X_2$ is a product of hyper-quadrics $Q_n\subset\mbb P^{n+1}$. The splitting criterion involves a very large number of cohomological conditions. By applying \ref{cor:relative}, we obtain:

\begin{m-corollary}\label{cor:pxq}
Let $X_1,X_2$ be as above. 
\begin{enumerate}
\item 
A vector bundle on $X_1$ splits if and only if it splits along a $\mbb P^2\times\ldots\times\mbb P^2\subset X_1$. (This reduces the number of cohomological tests to $3^t$.) 

\item 
A vector bundle on $X_1\times X_2$ splits if and only if it splits along some $X_1'\times X_2'\subset X_1\times X_2$, where $X_1'$ is a product of projective planes $\mbb P^2$ and $X_2'$ is a product of copies of $Q_3$. (The number of cohomological tests decreases dramatically.)
\end{enumerate}
\end{m-corollary}
\smallskip

\nit\textbf{(Vector bundles on products of minuscule varieties).} 
The minuscule homogeneous varieties are the following: the projective spaces, Grassmannians, spinor varieties, quadrics, the Cayley plane, and Freudenthal's variety. In \cite{hal2} it is proved that a vector bundle on a minuscule homogeneous variety $M$, $\dim M\ges 2$, splits if and only if its restriction to the union $M_2\subset M$ of the two-dimensional Schubert subvarieties splits. (It turns out that $M_2$ is either $\mbb P^2$ or a union of two copies of $\mbb P^2$ glued along a $\mbb P^1$.) The proof, which does not fit within the frame of this article, exploits the compatible F-splitting of the Schubert varieties and the properties of the minuscule weights. 

Rather \emph{surprisingly}, our approach reduces the problem concerning the splitting of vector bundles on products of minuscule varieties to the problem of splitting on products of projective $2$-planes. To our knowledge, there are no results in this direction, even for vector bundles on products of Grassmannians. Based on the articles \cite{ottv,mals}, one may speculate that, in the latter case, a cohomological splitting criterion should include a large number of tests, indexed by the Schur powers of the universal quotient bundles on the factors.  

\begin{m-theorem}\label{thm:minuscule}
Let $M^{(j)}$, $j=1,\dots,t$, be minuscule homogeneous varieties, $\dim M^{(j)}\ges 2$, and 
$X:=M^{(1)}\times\dots\times M^{(t)}.$ 

A vector bundle $\crl V$ on $X$ splits if and only if its restriction to $M^{(1)}_2\times\dots\times M^{(t)}_2$ splits, where each $M^{(j)}_2\subset M^{(j)}$ stands for the union of the $2$-dimensional Schubert subvarieties. 
\end{m-theorem}

\begin {m-proof}
The proof is by induction on $t+\rk\crl V$. For $t=1$, see \cite{hal2}. Let us prove the statement for $\tld X=X\times M$, with $X$ as above and $M$ minuscule: we assume that $\crl V_{X_2\times M_2}$ splits, where $X_2:=M^{(1)}_2\times\dots\times M^{(t)}_2\subset X$. The proof consists of two steps. 
\smallskip

\nit\textit{Claim 1}\quad $\crl V_{X\times M_2}$ splits. 
We prove this step by induction on $\rk(\crl V)$. 

Recall that $\mbb Z^{t+1}\cong\Pic(\tld X)\srel{\cong}{\to}\Pic(X_2\times M_2)$, so the line bundles on $\tld X$ are of the form 
$
\eO_{\tld X}(\unl{\tld\alpha})=
\eO_{M^{(1)}}(\alpha_1)\boxtimes\dots\boxtimes\eO_{M^{(t)}}(\alpha_t)\boxtimes\eO_M(k)=
\eO_X(\unl\alpha)\boxtimes\eO_M(k).
$ 
We deduce that 
\begin{equation}\label{eq:xm2}
\crl V_{X_2\times M_2}\cong
\uset{(\unl\alpha,k)\in\mbb Z^{t+1}}{\bigoplus} \eO_X(\unl\alpha)\boxtimes\eO_{M_2}(k)^{d_{\unl\alpha,k}}.
\end{equation}
Let us consider the diagram 
$$
\xymatrix@R=1.75em@C=3em{
X_2\times M_2\,\ar@{^(->}[r]\ar[d]^-f&X\times M_2\ar[d]^-f
\\ 
M_2\ar@{=}[r]&M_2.
}
$$
The induction hypothesis implies that $\crl V$ splits on the fibres of $f$, so
$$
\crl V_{X\times\{z\}}\cong\uset{\unl\alpha\in\mbb Z^t}{\bigoplus} 
\eO_X(\unl\alpha)^{d_{\unl\alpha}(z)},\; \forall\,z\in M_2.
$$ 
Actually, the multiplicities $d_{\unl\alpha}(z)$ are independent of $z\in M_2$. By restricting to $X_2\times\{z\}$ and by using \eqref{eq:xm2}, we deduce: 
$$
d_{\unl\alpha}:=\uset{k\in\mbb Z}{\sum}\;d_{\unl\alpha,k}
=d_{\unl\alpha}(z),\;\forall\,z\in M_2.
$$
Let $\unl a\in\mbb Z^t$ be a maximal element of $\{\unl\alpha\in\mbb Z^t\mid d_{\unl\alpha}\neq0\}$, for the lexicographic order. Then $f_*\big(\eO_X(-\unl a)\otimes\crl V\big)$ is locally free on $M_2$ of rank $d_{\unl a}$ and \eqref{eq:xm2} yields: 
$$
\crl R:=f_*\big(\eO_X(-\unl a)\otimes\crl V\big)\cong
\uset{k\in\mbb Z}{\bigoplus}\eO_{M_2}(k)^{d_{\unl a,k}}.
$$
(For this, observe that 
$H^0\big(X,\eO_X(\unl\alpha-\unl a)\big), H^0\big(X_2,\eO_X(\unl\alpha-\unl a)\big)\neq 0$ 
if and only if all the components of $\unl\alpha-\unl a$ are positive; the only such $\unl\alpha$ is $\unl a$ itself.) It follows that we have the exact sequence 
$\;0\to\eO_X(\unl a)\boxtimes\crl R\to\crl V\to\crl V'\to0.$

The first arrow is pointwise injective, so the quotient $\crl V'$ is locally free on $X\times M_2$, its restriction to $X_2\times M_2$ splits, and $\rk(\crl V')<\rk(\crl V)$. Hence, by the induction hypothesis, $\crl V'$ splits. Finally, we deduce that 
$$
\crl V\in\Ext^1(\crl V',\eO_X(\unl a)\boxtimes\crl R)=H^1\big(X\times M_2,\cHom(\crl V',\eO_X(\unl a)\boxtimes\crl R)\big)=0,
$$
because $X\times M_2$ is $1$-split and $\cHom(\crl V',\eO_X(\unl a)\boxtimes\crl R)$ is a direct sum of line bundles. We conclude by recurrence that $\crl V_{X\times M_2}$ splits. 
\smallskip


\nit\textit{Claim 2}\quad  $\crl V_{X\times M}$ splits.  

We denote by $M_d$ the union of all the $d$-dimensional Schubert subvarieties $S_d\subset M$. For any $(d+1)$-dimensional Schubert variety $S_{d+1}\subset M$, the intersection $\del S_{d+1}:=M_d\cap S_{d+1}$ is reduced and it is the union of the $d$-dimensional Schubert subvarieties of $S_{d+1}$; usually, it's called the boundary of $S_{d+1}$. 
With this notation, for $d\ges2$, the following properties hold (cf. \cite{hal2} and the references therein):  
\begin{itemize}[leftmargin=*]
\item 
$\mbb Z\cdot\eO_M(1)=\Pic(M)\to\Pic(S_d)$ is an isomorphism;
\item 
$\eO_M(1)\otimes \eO_{S_{d+1}}=\eO_{S_{d+1}}(\del S_{d+1});$
\item 
$S_d, \del S_{d+1}$ are $1$-split. 
\end{itemize}
Thus, for any Schubert variety $S_{d+1}\subset M$, the divisor $X\times\del S_{d+1}\subset X\times S_{d+1}$ is $1$-split and $\dim X$-positive. We deduce the implications: 
$$
\crl V_{X\times M_d}\mbox{ splits }
\;\srel{\del S_{d+1}\subset M_d}{\Longrightarrow}\;
\crl V_{X\times\del S_{d+1}}\mbox{ splits }
\;\ouset{\ref{thm:1}}{\text{Thm.}}{\Longrightarrow}\;
\crl V_{X\times S_{d+1}}\mbox{ splits.}
$$ 
Clearly, for $S'_{d+1}, S''_{d+1}\subset M$, the splittings of 
$\crl V_{X\times S'_{d+1}},\crl V_{X\times S''_{d+1}}$ coincide along the (reduced) intersection $X\times (S'_{d+1}\cap S''_{d+1})\subset X\times M_d$, so one gets a splitting of $\crl V_{X\times M_{d+1}}$. By repeating the argument, we deduce that $\crl V_{X\times M}$ splits.
\end {m-proof}

\begin{m-example}\label{expl:mr+bm}
Let $X=X^{(1)}\times\dots\times X^{(t)}$ be a product of Fano varieties of dimension at least four, with $\Pic(X^{(j)})\cong\mbb Z$ for all $j$. By applying the Theorem \ref{thm:1}, one reduces the problem of splitting of a vector bundle on $X$ to $S^{(1)}_3\times\dots\times S^{(t)}_3$, where each $S^{(j)}_3\subset X^{(j)}$ is an irreducible, $3$-dimensional complete intersection. 
\end{m-example}


\section{Splitting along divisors: a `probabilistic' approach}
\label{sct:global-gener}

In this section we obtain splitting criteria for vector bundles by restricting them to zero loci of \textit{generic} sections of \textit{globally generated}, partially positive line bundles. The global generation allows to replace the {\rm($2$-split)} by the weaker {\rm($1$-split)} condition, but we will have to consider very general test divisors instead of arbitrary ones. This explains the `probabilistic' attribute used in the title.

Note that, if $\eL$ is a $q$-ample line bundle on $X$ such that $\eL^d$ is globally generated for some $d\ges 1$, then $\eL$ is $q$-positive and the fibres of the morphism $f:X\to|\eL^d|$ are at most $q$-dimensional (cf. \cite[Theorem 1.4]{mats}), thus $\kappa(\eL)\ges\dim({\rm Image}(f))\ges\dim X-q$. Henceforth we replace $\eL^d$ by $\eL$.

Let the situation be as above. We start with general considerations: the equations defining $X,\eL,\crl V$ involve finitely many coefficients in $\kk$. By adjoining them to $\mbb Q$, we obtain a field extension $\mbb Q\hra\Bbbk$ of finite type (which depends on $\crl V$), so we may assume that $\Bbbk\hra\kk$; its algebraic closure is a \emph{countable}. After replacing $\Bbbk$ by $\bar{\Bbbk}$, we may assume that $X,\eL,\crl V$ are defined over a countable, algebraically closed subfield $\Bbbk$ of $\kk$; we denote these objects by $X_\Bbbk,\eL_\Bbbk,\crl V_\Bbbk$. 

The sheaf $\eG{:=}\,\Ker\big( H^0(X,\eL)\otimes\eO_X {\to}\,\eL\big)$ is locally free and the incidence variety 
$$
\cal D{:=}\{([s\,],x)\mid s(x)=0\}\subset|\eL|\times X
$$ 
is naturally isomorphic to the projective bundle $\mbb P(\eG)$ over $X$. The projections of $\cal D$ onto $|\eL|, X$ are denoted by $\pi,\rho$:
$$
\xymatrix@R=.75em{
&\;\cal D\;\ar[dl]^-\pi\ar[dr]_-\rho\ar@{^(->}[rr]&&|\eL|\times X
\\ 
|\eL|&&\,X.
}
$$
For any open $S\subset|\eL|$, let $\cal D_S:=\pi^{-1}(S)$; for $s\in|\eL|$, let $D_s:=\dvs(s)\subset X$.

All these objects are defined over $\Bbbk$, and they are denoted by $\eL_\Bbbk, \cal D_\Bbbk, \pi_\Bbbk, \rho_\Bbbk$. Let $K_{\kk}:=\kk(|\eL|)$ and $K_{\Bbbk}:=\Bbbk(|\eL_\Bbbk|)$ be the function fields of the projective spaces $|\eL|$, $|\eL_\Bbbk|$, respectively. 

\begin{m-definition}\label{def:very-gen}
We say that a property holds for a \emph{very general} point of some parameter space, if it holds on the complement of countably many proper subvarieties of that parameter space. In our case, we are interested in the splitting of $\crl V_{D_s}$, for $s\in|\eL|$ very general.
\end{m-definition}

\begin{m-lemma}\label{lm:very-gen} 
\begin{enumerate}
\item 
If $\crl V_{D_s}$ splits for a very general $s\in|\eL|$, then $(\rho^*\crl V)\otimes\bar{K}_\kk$ splits. 
\item If $(\rho^*\crl V)\otimes\bar{K}_\kk$ splits, then there is an \emph{analytic} open ball $\ball\subset|\eL|$ such that $D_s$ is smooth, for all $s\in\ball$, and ${(\rho^*\crl V)}_\ball$ splits over $\cal D_\ball$.
\end{enumerate} 
\end{m-lemma}

\begin {m-proof}
(i) Let $\tau:|\eL|\to |\eL_\Bbbk|$ be the trace morphism. Since $\Bbbk$ is countable and $\kk$ is not, $\tau(s)$ is the generic point of $|\eL_\Bbbk|$, for $s\in|\eL|$ very general. Hence $\Bbbk(\tau(s))=\Bbbk(|\eL_\Bbbk|)=K_\Bbbk$ and we obtain the Cartesian diagram: 
$$
\xymatrix@R=1.75em@C=3em{
D_s\ar[r]\ar[d]
&
\cal D_{\bar{K}_\Bbbk}
\ar[d]\ar[r]
&
\cal D_{\Bbbk}\ar[d]
\\  
\Spec(\kk)\ar[r]
&
\Spec(\bar{K}_\Bbbk)\ar[r]
&
\,|\eL_\Bbbk|.
}
$$ 
For varieties defined over algebraically closed fields ($\kk$ and $\bar{K}_\Bbbk$ in our case), the property of a vector bundle to be split commutes with base change. Then $\crl V_{\bar{K}_\Bbbk}$ splits on $\cal D_{\bar{K}_\Bbbk}$, so the same holds for $\crl V_{\bar{K}_\kk}$. 

\nit{(ii)} 
We note that $\crl V_{\cal D_{\bar{K}_\kk}}$ actually splits over an intermediate field $K_\kk\hra K'\hra\bar{K}_\kk$ finitely generated (and also algebraic) over $K_\kk$. Thus there is an open affine $S\subset|\eL|$, an affine variety $S'$ over $\kk$, and a \emph{finite morphism} $S'\srel{\si}{\to} S$ such that the direct summands of $\crl V_{\cal D_{\bar{K}_\kk}}$ are defined over $\kk[S']$; thus ${(\rho^*\crl V)}_{S'}$ splits. For $S$ sufficiently small, Bertini's theorem implies that $D_s$ is smooth, for all $s\in S$. Finally, there are open balls $\ball'\subset S'$ and $\ball\subset S$ such that $\si:\ball'\to\ball$ is an {\em analytic} isomorphism. Then the splitting of ${(\rho^*\crl V)}_{\ball'}$ descends to ${(\rho^*\crl V)}_{\ball}$ on $\cal D_{\ball}$.
\end {m-proof}

\begin{m-remark}\label{rmk:very-gen}
The previous lemma precises the meaning of a very general point $s\in|\eL|$: its coordinates should be algebraically independent over the definition field of $X,\eL,\crl V$. In particular, the notion of a very general point depends on $\crl V$ itself, actually on its field of definition. 

Often one wishes to have statements which involve \emph{general} points, rather than very general ones. The splitting of $\crl V_{D_s}$ for general $s\in|\eL|$ means, by definition, that $(\rho^*\crl V)\otimes K_\kk$ splits. This condition is \emph{more restrictive} than the splitting of $(\rho^*\crl V)\otimes \bar K_\kk$. 
\end{m-remark}

\begin{m-theorem}\label{thm:2} 
Suppose that $\eL\in\Pic(X)$ is globally generated, $(\dim X-4)$-positive, and $D\in|\eL|$ is very general (thus smooth).  If $X$ is $1$-split, then holds: $[\,\crl V$ splits $\Leftrightarrow$ $\crl V_D$ splits$\,]$.
\end{m-theorem}

The interest in this result is that it allows to test the splitting of vector bundles along divisors which are not $1$-split; this situation arises especially in low dimensions. Otherwise, of course, one applies the `deterministic' Theorem \ref{thm:1}. 

\begin {m-proof} 
Let $\ball$ be as in the Lemma \ref{lm:very-gen}. By \cite[Proposition 5.1]{ottem}, the cohomological dimension $\cd(X\sm D_s)\les \dim X-4$, for all $s\in\ball$, which implies that, for all $o,s,t\in\ball$, the intersections $D_{st}:=D_s\cap D_t$ and $D_{ost}:=D_o\cap D_s\cap D_t$ are non-empty and connected (cf. \cite[Ch. III, Corollary 3.9]{hart-as}). Note that the intersections are generically transverse, because $\eL$ is globally generated. 
\smallskip

\nit\textit{Claim~1}\quad 
Let $s,t\in\ball$ such that $D_s,D_t$ intersect transversally. Then holds: 
$$
\begin{array}{cl}
\res^X_{D_s}:\Pic(X)\to\Pic(D_s)&\mbox{is an isomorphism};
\\[1ex]
\res^{D_s}_{D_{st}}:\Pic(D_s)\to\Pic(D_{st})&\mbox{is injective}.
\end{array}
$$
The first statement is proved in the Theorem \ref{thm:pic-iso} and the second in the Proposition \ref{prop:pic-inj}. 
\smallskip

\nit\textit{Claim~2}\quad 
$\rho^*:\Pic(X)\to\Pic(\cal D_\ball)$ is an isomorphism. 
Indeed, fix $o\in\ball$ and consider the diagram: 
$$\xymatrix@R=.75em@C=4em{
&\Pic(D_s)\ar[dr]^-{\res^{D_s}_{D_{os}}}&
\\ 
\Pic(X)\ar[ur]^-{\res^X_{D_s}}
\ar[dr]_-{\res^X_{D_o}}\ar[rr]|{\;\res^X_{D_{os}}\;}
&&
\Pic(D_{os}).
\\ 
&\Pic(D_o)\ar[ur]_-{\res^{D_o}_{D_{os}}}&
}
$$
The composition $\Pic(X)\,{\srel{\;\rho^*}{\to}}\Pic(\cal D_\ball){\srel{\res_{D_o}}{\to}}\Pic(D_o)$ is bijective, so $\rho^*$ is injective. For the surjectivity, take $\ell\,{\in}\Pic(\cal D_\ball)$. If $\ell_{D_o}{\cong}\,\eO_{D_o}$, then $\ell_{D_s}\in\Pic^0(D_s)$, for all $s\in\ball$, so 
$$\{s\in\ball\mid\ell_{D_s}\not\cong\eO_{D_s}\}=\{s\in\ball\mid h^0(\ell_{D_s})=0\}$$ 
is open in $\ball$, hence 
$S:=\{s\in\ball\mid\ell_{D_s}\cong\eO_{D_s}\}$ 
is closed. 

On the other hand, the diagram above implies, by taking the restrictions to $D_{os}$, that $S$ contains all $s$ such that $D_s$ intersects $D_o$ transversally; thus $S$ is dense, so $S=\ball$. It follows that $\ell\cong\eO_{\cal D_\ball}$. For an arbitrary $\ell\in\Pic(\cal D_\ball)$, let $\eL\in\Pic(X)$ such that $\ell_{D_o}\cong\eL_{D_o}$, so ${\big((\rho^*\eL^{-1})\otimes\ell\big)}_{D_o}$ is trivial. This concludes the proof of the Claim~2.
\smallskip

Since $\crl V_\ball$ splits, we deduce that 
${(\rho^*\crl V)}_\ball\,{\cong}\,\rho^*\big(\underset{j\in J}{\bigoplus}\eL_j^{\oplus d_j}\big)$, 
with $\eL_j\in\Pic(X)$ pairwise non-isomor\-phic. We consider the following partial order on line bundles: 
$$
\eL\prec\euf M\;\Leftrightarrow\;\eL\not\cong\euf M 
\text{ and } H^0(\eL^{-1}\otimes\euf M)\neq0.
$$
For $s\in\ball$, let $J_{s,\mx}\subset J$ be the set of maximal elements for $\prec$ on $\Pic(D_s)\cong\Pic(X)$. By semi-continuity, the set $\{t\in\ball\mid J_{s,\mx}\subset J_{t,\mx}\}$ is open. Thus, after shrinking $\ball$, we may assume that $J_{s,\mx}\subset J$ is independent of $s$; we denote it by $J_\mx$.  

The maximality property implies that there is a natural, pointwise injective homomorphism
\begin{equation}\label{eq:h-inj}
h:\underset{\mu\in J_\mx}{\bigoplus}
\rho^*\eL_\mu\otimes
\underbrace{\pi^*\pi_*\big(\rho^*(\eL_\mu^{-1}\otimes\crl V)\big)} %
_{\srel{\circledast}{\cong} \;\; \eO_{\cal D_\ball}^{\oplus d_\mu}}
\to{\big(\rho^*\crl V\big)}_{\cal D_\ball}.
\end{equation}

\nit\textit{Claim~3}\quad \begin{minipage}[t]{.8\textwidth} 
$h$ descends to $X$, after a suitable base change in $\eO_{\cal D_\ball}^{\oplus d_\mu}$ by an analytic map 
$\beta:\ball\to\kern-1ex\underset{\mu\in J_\mx}{\prod}\kern-1.5ex\Gl(d_\mu)$.\end{minipage}

Indeed, we fix $o\in\ball$ and, for each $\mu\in J_\mx$, we fix a basis in 
$H^0(D_o,\eL_\mu^{-1}\crl V)\cong\kk^{d_\mu}$. (Bases are represented as square matrices whose columns are the vectors of the basis.) For any $s\in\ball$, $D_{os}$ is non-empty and connected, so there is a \emph{unique} basis in $H^0(D_s,\eL_\mu^{-1}\crl V)\cong\kk^{d_\mu}$ whose restriction to $D_{os}$ coincides with the restriction of the basis along $D_o$. (As $s\in\ball$ varies, the transition matrices from the trivialization $\circledast$ in \eqref{eq:h-inj} to these new bases yields the map $\beta$.)

We observe that, after this reparameterization, $h$ descends to the open set $\cU:=\rho(\cal D_{\ball})\subset X$. Indeed, define 
\begin{equation}\label{eq:bar-h}
\bar h:
\big(\underset{m\in J_\mx}{\mbox{$\bigoplus$}}\kern-1ex\eL_\mu^{\oplus d_\mu}\big)\otimes\eO_\cU
\to\crl V_\cU,\;\bar h(x):=h_s(x)\text{ for some }s\in\ball\text{ such that  }x\in D_s.
\end{equation}
In order to prove that $\bar h(x)$ is independent of $s\in\ball$, we must show that, for any $s,t\in\ball$, the restrictions to $D_{st}$ of the new bases in $H^0(D_s,\eL_\mu^{-1}\crl V),H^0(D_t,\eL_\mu^{-1}\crl V)$ coincide: it is enough to check this on the triple intersection $D_{ost}=D_o\cap D_{st}$ (which is non-empty, connected), where both bases are induced from $D_o$.

The homomorphism \eqref{eq:bar-h} yields the extension of locally free sheaves on $\cU$:
$$
0\to
\big(\underset{m\in J_\mx}{\bigoplus}\kern-1ex\eL_\mu^{\oplus d_\mu}\big)\otimes\eO_\cU
\to\crl V_\cU\to\crl W_\cU\to0,
\quad 
\rho^*(\crl W_\cU)\cong\kern-.5ex\rho^*\Big(
\underset{j\in J\sm J_\mx}{\bigoplus}\kern-1.5ex\eL_\mu^{\oplus d_\mu}
\otimes\eO_\cU\Big).
$$
The homomorphism on the left is pointwise injective. Recursively, we deduce that $\crl V_\cU$ is obtained as a successive extension of the line bundles $\eL_j\otimes\eO_\cU$, $j\in J$. (Note that, in the gluing process, we did not use that $\crl V$ is defined on all $X$; we used only its restriction to $\cU$.)

Since $\cU$ is an analytic neighbourhood of $D_o$, one gets induced extensions on the thickenings ${(D_{o})}_m$, $m\ges 0$, (cf. Definition \ref{def:D-thick}). But $X$ is $1$-split and $D_o$ is $(\dim X-4)$-ample, so 
\begin{equation}\label{eq:ext1}
0=\Ext^1(\eL,\euf M)\to \Ext^1\big(  \eL_{{(D_{o})}_m}, \euf M_{{(D_{o})}_m}\big)
\end{equation}
is an isomorphism, for all $\eL,\euf M\in\Pic(X)$, $m\gg0$. It follows that $\crl V_{{(D_o)}_m}$ splits, for $m\gg0$. By applying the Lemma \ref{lm:q-split}, we deduce that $\crl V$ splits on $X$. (Here it is necessary to have $\crl V$ defined on $X$.)
\end {m-proof}

\begin{m-example}\label{expl:p2p2} 
We have seen in \ref{expl:proj-bdl} that the splitting of a vector bundle on a projective bundle over a $1$-split variety can be verified along an arbitrary $\mbb P^2$-sub-bundle. At this stage, the reduction process given by the Theorem \ref{thm:1} stops. In this example we show that the Theorem \ref{thm:2} allows to decrease further the dimension of the test subvarieties. 

Let $(S,\eA)$ be an arithmetically Cohen-Macaulay surface, with $\eA$ ample and $\Pic(S)=\mbb Z\eA$. (A necessary and sufficient cohomological condition for the splitting of vector bundles on such surfaces has been obtained in \cite{ha+ta}.) The four-fold 
$$
X:=\mbb P(\eO_S\oplus\eA^{-m}\oplus\eA^{-m-n})\srel{f}{\to} S,\;\;m,n\ges0,
$$ 
is $1$-split. The line bundle $\eL:=\eO_f(1)\otimes f^*\eA$ is ample on $X$, and the general $D\in|\eL|$ is a smooth $\mbb P^1$-fibre bundle over $S$. In particular, $D$ is not $1$-split, so the theorem \ref{thm:1} can not be applied. However, the `probabilistic' Theorem \ref{thm:2} still applies: a vector bundle $\crl V$ on $X$ splits if and only if it does split on a very general $D$. 
\end{m-example}


\section{Triviality criteria}\label{sct:triv}

Finally, in this section we restrict our discussion to the case of the trivializable vector bundles. The motivation is, first, that the general `effective splitting criterion' \ref{prop:h1} is not explicit enough, especially for partially ample line bundles which are pulls-back (cf. Remark \ref{rmk:c}, Theorem \ref{thm:rel-ample}). Second, it is desirable to remove the $1$-, $2$-split conditions which appear throughout the sections \ref{sct:q-ample}, \ref{sct:global-gener}, and which are imposed precisely to ensure the vanishing \eqref{eq:h1}. 

Unfortunately, the Kodaira vanishing does not hold for $q$-ample line bundles. Thus, to obtain effective results in this situation, one must find appropriate conditions which imply the Kodaira vanishing. These lines of thought lead to the triviality criteria below.


\begin{m-lemma}\label{lm:pic-inj}
Assume that $\eL\in\Pic(X)$ is $(\dim X-2)$-ample and satisfies $H^i(X,\eL^{-a})=0$, for all $a\ges 1$ and $i=0,1,2$. For any vector bundle $\crl V$ on $X$ and $D\in|\eL|$, one has the equivalence: 
$\;[\,\crl V\cong\eO_X^{\oplus r}\;\Leftrightarrow\;\crl V_D\cong\eO_D^{\oplus r}\,].$
\end{m-lemma}

\begin {m-proof}
The hypothesis implies that $H^i(\eL_D^{-a})=0$, for all $a\ges 1$, $i=0,1$, and  $H^0(\eO_{D_a})=\kk$, for $a\ges0$, so $H^0(\crl E)\to H^0(\crl E_D)=\End(\kk^r)$ is an isomorphism, by the Proposition \ref{prop:h1}. Hence $\crl V$ splits, actually $\crl V\cong\euf M^{\oplus r}$ for some $\euf M\in\Pic(X)$. As $\crl V_D=\eO_D^{\oplus r}$, we deduce that both $\euf M_D$ and $\euf M_D^{-1}$ admit non-trivial sections, so $\euf M_D\cong\eO_D$. According to the Proposition \ref{prop:pic-inj}, $\Pic(X)\to\Pic(D)$ is injective, so $\euf M$ is trivial. 
\end {m-proof}


\begin{m-theorem}\label{thm:q-global-gener}
Consider a vector bundle $\crl V$ on $X$, $\eL\in\Pic(X)$, and $D\in|\eL|$. In any of the following cases, we have:
$\;[\,\crl V\cong\eO_X^{\oplus r}\;\Leftrightarrow\;\crl V_D\cong\eO_D^{\oplus r}\,].$ 
\begin{enumerate}
\item[(a)] 
$\eL\in\Pic(X)$ is semi-ample and $(\dim X-3)$-ample; 
\item[(b)] 
$\eL$ is relatively ample for a morphism $X\srel{f}{\to}Y$, with $\dim(X)-\dim(Y)\ges3$. 
\end{enumerate}
\end{m-theorem}

\begin {m-proof}
In both cases, the Theorem \ref{thm:global-gen} implies $H^i(X,\eL^{-a})=0$ for $a\ges 1$, $i=0,1,2$, so we can apply the Lemma \ref{lm:pic-inj}.
\end {m-proof}


\subsection{The case of Frobenius split (F-split) varieties}\label{ssct:frob} 
These objects are ubiquitous, especially in representation theory. Examples of F-split varieties (defined in characteristic zero) include Fano varieties (cf. \cite[Exercise 1.6E(5)]{brion-kumar}), spherical varieties, in particular projective homogeneous varieties and toric varieties (cf. \cite{br-in}, \cite[Section 31]{tim}). The notions and properties which are relevant for us are summarized in the appendix \ref{sct:frob}. 

\begin{m-theorem}\label{thm:F-q-split}
Let $D$ be a $(\dim X-3)$-ample, effective divisor, which is F-split. Then holds: 
$$
\crl V\cong\eO_X^{\oplus r}
\;\Leftrightarrow\;
\crl V_D\cong\eO_D^{\oplus r}.
$$  
\end{m-theorem}

The F-splitting allows to handle more `exotic' situations. Many examples arise from varieties $X$ which are compatibly split with respect to a divisor $D$. 

\begin {m-proof}
Since $\eO_D(D)$ is $(\dim D-2)$-ample, the Theorem \ref{thm:F-split} 
implies $H^1(D,\eO_D(-aD))=0$, for all $a\ges 1$. The conclusion follows from the Proposition \ref{prop:h1}. 
\end {m-proof}

\begin{m-corollary}\label{cor:f-split}
Let $X$ be a smooth projective variety whose anti-canonical line bundle $\omega_X^{-1}$ is $(\dim X-3)$-ample. Assume that $X$ is F-split by $\si\in H^0(\omega_X^{-1})$, and denote $D:=\dvs(\si)$. Then holds: 
$[\,\crl V\cong\eO_X^{\oplus r}\;\Leftrightarrow\;\crl V_D\cong\eO_D^{\oplus r}\,].$
\\ The criterion applies, in particular, in the following cases: 
\begin{enumerate}
\item[(a)] $X$ is a Fano variety of dimension at least three;
\item[(b)] $X$ is a spherical variety whose anti-canonical bundle is $(\dim X-3)$-ample. 
\end{enumerate}
\end{m-corollary}

\begin {m-proof}
The hypothesis implies that $D$ is F-split, compatibly with the splitting defined by $\si$. 

Fano and spherical varieties are Frobenius split, compatibly with suitable anti-canonical divisors (cf. \cite[1.6.E(5), pp. 56]{brion-kumar} and \cite[Theorem 1]{br-in}, respectively). 
\end {m-proof}

\subsection{The case of toric varieties}\label{toric} 
A non-trivial application of the ideas developed inhere arises when $X:=X_\Si$ is the smooth projective toric variety defined by the regular fan $\Si$. 

\begin{m-remark}\label{rmk:toric}
{\rm(i)} 
$X$ is F-split, compatibly with the invariant divisors $D_\rho, \rho\in\Si(1)$, and their intersections (cf. \cite[Exercise 1.3E(6)]{brion-kumar}). 
\\{\rm(ii)} Let 
\begin{equation}\label{eq:dta}
\Delta:=\underset{\rho\in\Si(1)}{\mbox{$\sum$}}{D_\rho}
\end{equation}
be the torus-invariant, anti-canonical divisor. Its complement is $X\sm\Delta\cong(\mbb C ^*)^{\dim X}$, so the cohomological dimension equals $\cd(X\sm\Delta)=0$. 
\end{m-remark}

\begin{m-theorem}\label{thm:split-toric} 
Let $\crl V$ be an arbitrary vector bundle on the smooth toric variety $X$. The following statements hold: 
\begin{enumerate}
\item 
Let $\disp\hat X:=\varprojlim_m\Delta_m$ be the formal completion of $X$ along $\Delta$. If $\dim X\ges 2$, then one has the equivalence: 
$
\;[\,\crl V\mbox{ splits}\;\Leftrightarrow\;\crl V\otimes\eO_{\hat X}\mbox{ splits}\,].
$
\item 
Assume that 
$\dim X,\;\codim\big(\bloc(\omega_{X}^{-1})\,\big) 
\ges 3.$ 
Then we have: 
\begin{enumerate}
\item[(a)] 
$[\,\crl V\mbox{ splits}\;\Leftrightarrow\;\crl V_{\Delta_m}\mbox{ splits}\,]$, for $m\gg0$. 
(See \ref{rmk:c}(i) for a lower bound on $m$.)
\item[(b)] 
$[\,\crl V\mbox{ is trivial}\;\Leftrightarrow\;\crl V_{\Delta}\mbox{ is trivial}\,]$.
\end{enumerate}
\end{enumerate}
\end{m-theorem}

\begin{m-proof}
\nit(i) See the Proposition \ref{prop:formal}. 

\nit(ii) The Theorem \ref{thm:q-ample} implies that $\omega_{X}^{-1}$ is $(\dim X-3)$-ample. Now the conclusion follows from the Proposition \ref{prop:h1} and the Corollary \ref{cor:f-split}, respectively. 
\end{m-proof}

\begin{m-remark}
\nit{\rm(i)} 
One may wonder if it is possible to have a splitting criterion for toric varieties which involves an \emph{irreducible}, torus-invariant, $(\dim X-2)$-ample divisor.  In general, the answer is ``no''; \emph{reducible} divisors are necessary for the following reason.  If $D$ is such an irreducible divisor, then ${\rm cd}(X\sm D)\les\dim X-2$, so $D$ must intersect all the other torus-invariant divisors. Hence $\Si$ has the following property: if $\xi_D\in\Si(1)$ defines $D$, then $\xi_D,\xi$ span a cone of $\Si$, for all $\xi\in\Si(1)\sm\{\xi_D\}$. This condition is clearly not satisfied in general.
\smallskip

\nit{\rm(ii)} 
It is rather surprising, but the issue concerning the bare existence of non-trivial vector bundles on toric varieties is not settled yet, in general (cf. \cite{gharib+karu}).
\end{m-remark}


\appendix

\section{About $q$-ample and $q$-positive line bundles} \label{sct:setup}

In this section we summarize the notions and the results about partial positivity for line bundles which are 
used in this note. Henceforth $X$ stands for a smooth projective variety of arbitrary dimension, defined over $\kk$. 

\begin{m-definition}\label{def:q-ample} 
Consider a line bundle $\eL$ on $X$. 

\begin{enumerate}
\item 
(cf. \cite{totaro}) $\eL$ is \emph{$q$-ample}, if for any coherent sheaf $\eF$ on $X$ there is $m_\eF>0$ such that $H^i(X,\eF\otimes\eL^m)=0$, for all $m\ges m_\eF$ and $i>q$. 

\item 
(cf. \cite{ag,dps}) $\eL$ is \emph{$q$-positive}, if it admits a Hermitian metric whose curvature is positive definite on a subspace of $T_{X,x}$ of dimension at least $\dim X-q$, for all $x\in X$;  equivalently, the curvature has at each point $x\in X$ at most $q$ negative or zero eigenvalues. 

\item 
$\eL$ is \emph{semi-ample}, if a tensor power of it is globally generated. 

\item 
Assume that $\exists\,a\ges1$ such that $H^0(X,\eL^a)\neq0$. The \emph{Kodaira-Iitaka dimension} and the \emph{(stable) base locus} of $\eL$ are defined as follows (cf. \cite[Section 2.1]{laz1}):
$$ 
\begin{array}{l}
\disp
\kappa(\eL):=\mbox{transcend.\,deg.}_{\kk}
\big(\,\uset{a\ges 0}{\bigoplus} H^0(X,\eL^a)\big){-}1
=\max_{a\ges 1} 
\dim\big({\rm Image}(X\dashto |\eL^a|)\big);
\\[2ex] 
\bloc(\eL):=\kern-2ex\uset{s\in H^0(X,\eL)}{\bigcap}\kern-2ex\dvs(s),\;\; \sbloc(\eL):=\uset{a\ges 1}{\bigcap}\bloc(\eL^a)_{\rm red}.
\end{array}
$$
\end{enumerate}
\end{m-definition}

\begin{m-remark}\label{rmk:tensor}
(i) Any $q$-positive line bundle is $q$-ample (cf. \cite[Proposition 28]{ag}, \cite[Proposition 2.1]{dps}), but the converse is false (cf. \cite[Theorem 10.3]{ottem}). 

\nit(ii) If $\ell,\eL\in\Pic(X)$ and $\eL$ is $q$-ample, then $\ell\otimes\eL^m$ is $q$-ample for $m\gg0$. This is a direct consequence of the uniform $q$-ampleness property \cite[Theorem 6.4]{totaro}. 

\nit(iii) The Definition \ref{def:q-ample}(i) makes sense for any projective scheme, not necessarily 
smooth. This more general situation occurs in \ref{thm:F-split}. 

\nit(iv) If $\eL$ is $q$-ample (positive), then it is also $q'$-ample (positive), for any $q'\ges q$; the larger the value of $q$ the weaker the restriction on $\eL$. \textit{E.g.} the $\dim X$-ampleness (positivity) is an empty condition.
\end{m-remark}

\begin{m-theorem}\label{thm:mats} 
Assume that $\eL$ is semi-ample and $q$-ample. Then $\eL$ is $q$-positive and 
$$\dim X\les q+\kappa(\eL).$$
\end{m-theorem}

\begin {m-proof}
See \cite[Theorem 1.4]{mats}. For suitable $a$, the image of the morphism $X\to|\eL^a|$ is $\kappa(\eL)$-dimensional; we denote the image by $Y$. Hence $\dim X-\kappa(\eL)$ equals the dimension of the generic fibre of $X\to Y$; by \textit{loc. cit.}, the dimension of all the fibres is bounded above by $q$. 
\end {m-proof}

\begin{m-lemma}\label{lm:q-rel}
Let $X,Y$ be smooth projective varieties, and $f:X\to Y$ be a smooth, surjective morphism of relative 
dimension $\delta$. Then the following implications hold:
\begin{enumerate}
\item 
If $\euf M\in\Pic(Y)$ is $q$-ample, then $\eL:=f^*\euf M$ is $(\delta+q)$-ample;
\item 
If $\euf M\in\Pic(Y)$ is $q$-positive, then $\eL:=f^*\euf M$ is $(\delta+q)$-positive;
\item 
If $\eL\in\Pic(X)$ is $f$-relatively ample, then $\eL$ is $\dim Y$-positive. 
\end{enumerate}
\end{m-lemma}

\begin {m-proof}
Leray's spectral sequence implies (i); for (ii), the pull-back metric on $\eL$ satisfies \ref{def:q-ample}.

\nit(iii) Note that $\eL':=\eL\otimes f^*\eA$ is ample, for $\eA\in\Pic(Y)$ sufficiently ample. Then $m\eL'$, $m\gg0$, defines an embedding $\iota:X\to \mbb P$ into some projective space; the morphism $(f,\iota):X\to Y\times\mbb P$ is an embedding too, and $\eL^m=(f,\iota)^*\big(\eA^{-m}\boxtimes\eO_{\mbb P}(1)\big)$. The restriction to $X$ of the product metric on $\eA^{-m}\boxtimes\eO_{\mbb P}(1)$ is positive definite on $\Ker(\rd f)$. 
\end {m-proof}

\begin{m-theorem}\label{thm:global-gen}
\begin{enumerate}
\item 
Assume that $\eL\in\Pic(X)$ is \emph{semi-ample} and \emph{$q$-ample}, $q\les\dim X-1$. Then holds $\;H^i(X,\eL^{-a})=0,\;\forall\,i\les\dim X-q-1,\;\forall\,a\ges 1.$ 
\item \emph{(relative Kodaira vanishing)} 
Consider a morphism $X\srel{f}{\to} Y$ with $Y$ projective, and let $\eL\in\Pic(X)$ be \emph{$f$-relatively ample}. Then holds $\;H^i(X,\eL^{-1})=0,\,\forall\,i<\dim X-\dim Y.$ 
\item 
Let $Z$ be a projective, equidimensional, reduced, weakly normal variety which satisfies the condition $S_2$ of Serre. Let $\eL\in\Pic(Z)$ be relatively ample for $Z\srel{f}{\to}Y$, with $Y$ projective, and $\dim Z-\dim Y\ges 2$. Then holds $H^1(Z,\eL^{-1})=0$. 
\end{enumerate}
\end{m-theorem}

\begin {m-proof}
(i) The Grauert-Riemenschneider theorem (cf. \cite[Theorem 7.73]{sh+so}) yields the vanishing of $H^i(X,\eL^{-a})$, for all  $a\ges 1$ and $i\les\kappa(\eL)-1$. Now, we use the inequality in \ref{thm:mats}. 

\nit(ii) The claim follows from \cite[Theorem 1-2-3]{kmm}: $R^if_*(K_X(D))=0,\;\forall i>0.$ 
Indeed, for $i<\dim X-\dim Y$, the Leray-spectral sequence implies that we have:
$$
H^i(X,\eO_X(-D))\cong H^{\dim X-i}(X,K_X(D))=H^{\dim X-i}(\,Y,f_*(K_X(D))\,)=0.
$$
An independent proof, can be obtained by following the lines of (iii) below. 

\nit(iii) The previous argument does not apply because $Z$ is not necessarily smooth. 
Without loss of generality, we may assume that $f$ is surjective. For $\dim Y=0$, the vanishing is proved in \cite[Theorem 3.1]{ar-ja}. Now take a very ample line bundle $\eA$ on $Y$, such that $\eL\otimes f^*\eA$ is ample on $Z$. Bertini's theorem implies that $Z':=f^{-1}(Y')$ satisfies the same assumptions as $Z$, for general $Y'\in|\eA|$ (cf. \cite[Theorem 1]{cgm2}, also \cite[Corollary 1.9]{cgm1}). Finally, observe that in the exact sequence  
$$
\dots\to H^1(Z,(\eL\otimes\eA)^{-1})\to H^1(Z,\eL^{-1})\to H^1(Z',(\eL\otimes\eO_{Z'})^{-1})\to\dots
$$
both extremities vanish: by \cite[Theorem 3.1]{ar-ja}, for the left-hand-side, and by the induction hypothesis, for the right-hand-side. 
\end {m-proof}

\begin{m-proposition}\label{prop:pic-inj}
Assume that $\eL\in\Pic(X)$ is $(\dim X-2)$-ample and satisfies 
$$
H^i(X,\eL^{-a})=0,\;\forall\,a\ges 1\mbox{ and }i=0,1,2.
$$ 
Then, for any $D\in|\eL|$, the restriction $\Pic(X)\to\Pic(D)$ is injective. 

The vanishing condition is satisfied if $\eL$ is semi-ample with $\kappa(\eL)\ges 3$, in particular for $\eL$ semi-ample and $(\dim X-3)$-ample. 
\end{m-proposition}

\begin {m-proof}
The exact sequence $0\to \eL^{-1}\to\eO_X\to\eO_{D}\to0$ implies $H^i(\eL_D^{-a})=0$, for all $a\ges 1$ and $i=0,1$. By plugging this into 
$0\to \eL_D^{-a}\to\eO_{D_a}^\times\to\eO_{D_{a-1}}^\times\to0,$
we deduce that $\Pic(D_a)\to\Pic(D_{a-1})$ is injective, for all $a\ges 1$. 

Take $\euf M\in\Ker(\Pic(X)\to\Pic(D))$, so $\euf M_D\cong\eO_D$; it follows $\euf M_{D_a}\cong\eO_{D_a}$, for all $a\ges 0$. Note that $\mbb C\cong H^0(\eO_X)\to H^0(\eO_{D_a})$ is an isomorphism, for all $a\ges 0$, so $H^0(\euf M_{D_a})\cong\kk\cong H^0(\euf M_{D_a}^{-1})$. But the restrictions $H^0(\euf M)\to H^0(\euf M_{D_a})$ and $H^0(\euf M^{-1})\to H^0(\euf M^{-1}_{D_a})$ are isomorphisms, for $a\gg0$, thus $\euf M\cong\eO_X$. The final statement follows from the Grauert-Riemenschneider theorem \cite[Theorem 7.73]{sh+so} and the inequality in \ref{thm:mats}. 
\end {m-proof}

\begin{m-theorem}\label{thm:pic-iso}
Let $D\subset X$ be an effective divisor. The restriction $\Pic(X)\to\Pic(D)$ is an isomorphism in the cases enumerated below: 
\begin{enumerate}
\item[(a)] 
$D$ is \emph{smooth} and $q$-positive, with $q\les\dim X-4$; 
\item[(b)] \emph{(relative Picard-Lefschetz)} 
$D$ is \emph{arbitrary}, relatively ample for a morphism $X\srel{f}{\to}Y$, with $Y$ projective, and $\dim X-\dim Y\ges 4$.
\end{enumerate}
\end{m-theorem}

\begin {m-proof}
(i) The $q$-positivity of $\eL$ implies that $H^i(X;\mbb Z)\to H^i(D;\mbb Z)$, $i\les2$, are isomorphisms (cf. \cite[Theorem III]{bott-lefschetz}, \cite[Lemma 10.1]{ottem}), so the same holds with $\kk$-coefficients. The Hodge decomposition for $X,D$ yields $H^i(X;\eO_X)\srel{\cong}{\to} H^i(D;\eO_D)$, for $i\les2$. By comparing the exponential sequences for $X,D$, we deduce that  $\Pic(X)\srel{\cong}{\to}\Pic(D)$. 

\nit(ii) We may assume that $f$ is surjective. For $\dim Y=0$, this is the Picard-Lefschetz theorem \cite[IV\S3, Theorem 3.1]{hart-as}. Now we make the inductive step. Let $\eA$ be a very ample line bundle on $Y$, such that $\eO_X(D)\otimes f^*\eA$ is ample on $X$. Bertini's theorem implies that the general hyperplane section $Y'\subset Y$ has the following properties: 
\begin{itemize}[leftmargin=*]
\item 
$X':=f^{-1}(Y')$ is smooth;
\item 
$X'$ does not contain the support of any irreducible component of $D$, that is the schematic intersection $D\cdot X'$ of $D,X'$ is a divisor in $X'$. 
\end{itemize}

\nit The divisor $D+X'$ is ample on $X$, so the Picard-Lefschetz theorem (cf. \textit{idem}) implies that $\Pic(X)\to\Pic(D+X')$ is an isomorphism. Since $X'$ intersects $D$ properly, a line bundle (an invertible sheaf) $\ell\in\Pic(D+X')$ is uniquely determined by: 
\begin{itemize}[leftmargin=*]
\item 
a pair $(\ell_D,\ell_{X'})\in\Pic(D)\times\Pic(X')$;
\item 
an isomorphism $\ell_D\otimes\eO_{D\cdot X'}\cong\ell_{X'}\otimes\eO_{D\cdot X'}$.%
\footnote{
Let $\cU\subset\Spec( \mbb C[\xi_1,\dots,\xi_{\dim X-1},y])$ be a local analytic chart on $X$, such that $\{y=0\}$ and $\{f(\unl{\xi},y)=0\}$ are the local equations of $X'$ and $D$, respectively. Then 
$\disp 
0\to\frac{ \mbb C[\underline{\xi},y] }{ \lran{y\cdot f} }
\to\frac{\mbb C[\unl{\xi},y]}{\lran{y}}\oplus\frac{\mbb C[\unl{\xi},y]}{\lran{f}}
\to\frac{\mbb C[\unl{\xi},y]}{\lran{y,f}}\to0
$ 
shows that the germs of regular functions on $D+X'={\rm Var}(y\cdot f)$ consist of pairs of regular functions on $D,X'$ which agree an $D\cdot X'={\rm Var}(y,f)$.
}
\end{itemize} 
By the induction hypothesis, $\Pic(X')\to\Pic(D\cdot X')$ is an isomorphism, hence 
$\Pic(D+X')\to\Pic(D)$ is an isomorphism too. The conclusion follows.
\end {m-proof}


The precise condition for the $(\dim X-1)$-ampleness of a line bundle on $X$ is given in \cite[Theorem 9.1]{totaro}. In this article, we often assume that some effective divisor $D\subset X$ is either $(\dim X-3)$ or $(\dim X-2)$-ample, so we need a practical method to verify this condition. Note that $\eO_X(D)$ is automatically $(\dim X-1)$-ample. 

\begin{m-theorem}\label{thm:q-ample}
Let $\Delta$ be an effective divisor on $X$ and  $\eL:=\eO_X(\Delta)$. We assume: 
$$
\begin{array}{rl}
\rm(i)&
\cd(X\sm\Delta)\les\dim X-c,
\\[1ex] 
\rm(ii)&
\codim\big(\sbloc(\eL)\,\big) 
\ges c,
\end{array}\quad \mbox{for }c\ges1.
$$
Then $\eL$ is $(\dim X-c)$-ample.
\end{m-theorem}
Here $\cd(\cdot)$ stands for the cohomological dimension. Recall that, according to \cite[Proposition 5.1]{ottem}, if $\Delta$ is $q$-ample then $\cd(X\sm\Delta)\les q$. 

\begin{m-proof}
Let us analyse the assumption (ii). The statement of the theorem is invariant after replacing $\eL$ by a positive power $\eL^a$, thus we may assume that 
$$
\sbloc(\eL)=\bloc(\eL)_{\rm red}.
$$ 
The codimension of the former is at least $c$, so there exist at least $c$ algebraically independent sections in $\eL$, that is $\kappa(\eL)\ges c-1$. Bertini's theorem implies that, for general divisors $D_1,\dots,D_{c-1}\in|\eL|$, the schematic intersection $Z_k:=D_1\cdot\ldots\cdot D_k,\;k\les c-1,$ has codimension $k$ in $X$. Furthermore, since $\dim(\bloc(\eL))\les\dim Z_{c-1}-1$, there is a section in $\eL$ which vanishes along a non-trivial divisor $Z_c\subset Z_{c-1}$; otherwise, some component of $Z_{c-1}$ would be contained in $\bloc(\eL)$. 

We deduce the following exact sequences of sheaves on $X$, for $k=1,\dots,c$: 
\begin{equation}\label{eq:Zk}
0\to 
\eL^{m-1}\otimes\eO_{Z_{k-1}}\to
\eL^{m}\otimes\eO_{Z_{k-1}}\to
\eL^{m}\otimes\eO_{Z_{k}}\to
0,
\quad\mbox{($Z_0:=X$)}.
\end{equation}
Now we use (i): since $\cd(X\sm\Delta)\les \dim X-c$, any coherent sheaf $\eG$ on $X$ satisfies (cf. \cite[Equation (5.1)]{ottem}): 
\begin{equation}\label{eq:Hi=0}
\varinjlim_m H^i(X,\eG\otimes\eL^m)=
H^i(X\sm\Delta,\eG)=0,\;\forall\,i>\dim X-c.
\end{equation}
The proof of the theorem is by induction on $c$. Recall that it is enough to check the $q$-ampleness property for locally free sheaves $\eF$ on $X$; we fix one. 
For $c=1$, we tensor by $\eF$ the sequence \eqref{eq:Zk}, $k=1$, and obtain 
$$
H^{\dim X}(X,\eF\otimes\eL^{m-1})\to H^{\dim X}(X,\eF\otimes\eL^m)\to0.
$$
Thus the sequence of the cohomology groups eventually becomes stationary. By inserting into \eqref{eq:Hi=0}, we deduce that $H^{\dim X}(X,\eF\otimes\eL^m)=0$, for $m\gg0$. 

Now we proceed with the inductive step: assume that the theorem holds for $c$ and prove it for $c+1$. So we must show that $H^{\dim X-c}(X,\eF\otimes\eL^m)=0$, for $m\gg0$. 
\begin{itemize}[leftmargin=*] 
\item[--] 
The sequence \eqref{eq:Zk}, $k=c+1$, tensored by $\eF$ yields: 
\\[1ex] 
$
H^{\dim X-c}\big(X,(\eF\otimes\eO_{Z_{c}})\otimes\eL^{m-1}\big)\to
H^{\dim X-c}\big(X,(\eF\otimes\eO_{Z_{c}})\otimes\eL^m\big)\to0,
$ 
\\[1ex] 
because $\dim Z_{c+1}=\dim X-c-1$. For $\eG=\eF\otimes\eO_{Z_{c}}$, \eqref{eq:Hi=0} implies, the same as before, that 
$H^{\dim X-c}\big(X,(\eF\otimes\eO_{Z_{c}})\otimes\eL^m\big)=0,$ for $m\gg0$. 
\item[--] 
Insert this into the long sequence in cohomology corresponding to \eqref{eq:Zk}, $k=c$, and find for $m\gg0$:
\\[1ex] 
$
H^{\dim X-c}\big(X,(\eF\otimes\eO_{Z_{c-1}})\otimes\eL^{m-1}\big)\to
H^{\dim X-c}\big(X,(\eF\otimes\eO_{Z_{c-1}})\otimes\eL^m\big)\to0.
$
\\[1ex] 
Then \eqref{eq:Hi=0}, with $\eG=\eF\otimes\eO_{Z_{c-1}}$, yields 
$H^{\dim X-c}\big(X,(\eF\otimes\eO_{Z_{c-1}})\otimes\eL^m\big)=0$,  for $m\gg0$. 
\item[--] 
Repeat this procedure---use \eqref{eq:Zk}, for $k=c-1,\dots,1$, and \eqref{eq:Hi=0}---until we get $$H^{\dim X-c}(X,\eF\otimes\eL^m)=0,\text{ for }m\gg0.$$
\end{itemize}
This completes the proof of the theorem.
\end{m-proof}


\section{About Frobenius split (F-split) varieties}\label{sct:frob}

We recall the relevant definitions; the reference for the concept of Frobenius splitting is the book \cite{brion-kumar}, and also \cite[Section 31]{tim} for applications. 

\begin{m-definition}(cf. \cite[Definition 1.1.3 and Section 1.6]{brion-kumar}) 
Let $Z_p$ be a projective variety over $\bar{\mbb F}_p$ (the algebraic closure of the field $\mbb Z/p\mbb Z$). The absolute Frobenius morphism $F$ of $Z_p$ determines the sheaf homomorphism $F^\sharp:\eO_{Z_p}\to F_*\eO_{Z_p}$. One says that $Z_p$ is \emph{F-split} if  there is an $\eO_{Z_p}$-linear homomorphism 
$$
\vphi:F_*\eO_{Z_p}\to\eO_{Z_p}\text{ such that }
\vphi\circ F^\sharp=\bone_{\eO_{Z_p}}.
$$
A closed subscheme $Y\subset X$ defined by the sheaf of ideals $\eI_Y$ is \textit{compatibly split}, if $\vphi(\eI_Y)=\eI_Y$. 
\end{m-definition}

If $Z$ is a projective variety defined over a field of characteristic zero, there is a finite set $\mathbf{s}$ of primes, a finitely generated $\mbb Z[\mathbf{s}^{-1}]$-algebra $R$, and a smooth $\Spec(R)$-scheme $\cal Z$ such that $Z=\cal Z\times_{R}\kk$. If $\eL\in\Pic(Z)$, one may choose $R$ in such a way that $\eL$ also extends over $\Spec(R)$. 

\begin{m-definition}\label{def:red-p}
For a maximal ideal $\mfrak m\in\Spec(R)$, the residue field $k(\mfrak m)$ is a finite extension of 
$\mbb F_p$, with $p\not\in\mathbf{s}$. 
The variety $Z_p:=\cal Z\times_R\ovl{k(\mfrak m)}$ is called \emph{a reduction modulo $p$} 
of $Z$. (Note that $\ovl{k(\mfrak m)}\cong\ovl{\mbb F}_p$)

We say that $Z$ is \emph{F-split} if $Z_p$ is so, at infinitely many $\mfrak m\in\Spec(R)$. (Such a subset is automatically dense in $\Spec(R)$.) 
\end{m-definition}

In our context, the importance of the Frobenius splitting is captured in the following Kodaira vanishing theorem for $q$-ample line bundles, which (apparently) has not been observed so far. 

\begin{m-theorem}\label{thm:F-split}
Let $Z$ be a projective, equidimensional, Cohen-Macaulay, F-split variety over $\kk$, and let $\eL\in\Pic(Z)$ be $q$-ample. Then holds $H^i(Z,\eL^{-1})=0$, for all $i<\dim Z-q$. 
\end{m-theorem}
In this note, the result will be applied in the case where $Z$ is a compatibly split, normal crossing divisor of a smooth variety $X$. 

\begin {m-proof}
Consider $\cal Z\srel{\pi}{\to}\Spec(R)$ as above, such that $\eL$ extends to $\crl L\to\cal Z$. Then for all primes $p$ large enough, $\crl L_p\in\Pic(Z_p)$ is still $q$-ample (cf. \cite[Theorem 8.1]{totaro}), so $H^i(Z_p,\crl L_p^{-m})=0$, for $i<\dim Y-q$ and $m\gg0$, by Serre duality (cf. \cite[Ch. III, Corollary 7.7]{hart}). The F-splitting property implies that $H^i(Z_p,\crl L_p^{-1})=0$ (cf. \cite[Lemma 1.2.7]{brion-kumar}). Finally, the generic rank of the coherent sheaf $R^i\pi_*\crl L^{-1}$ on $\Spec(R)$ is constant, and the conclusion follows from the fact that the vanishing holds for infinitely many primes $p$. 
\end {m-proof}

\begin{m-remark}\label{rmk:expl-F}
{\rm(i)} The F-splitting of a non-singular variety $X_p$ (defined in characteristic $p$) is given by an element in $H^0(X_p,\omega_{X_p}^{1-p})$ satisfying a certain algebraic equation, where $\omega_{X_p}$ stands for the canonical sheaf (cf. \cite[Theorem 1.3.8]{brion-kumar}). 
\\{\rm(ii)} 
In \emph{characteristic zero}, an important source of F-splittings arise from varieties $X$ which have the property that their reduction $X_p$ modulo $p$ is split by the $(p-1)$-st power of (the mod $p$ reduction of) a section $\si\in H^0(X,\omega_X^{-1})$; in this case $D:=\dvs(\si)$ is a compatibly split subvariety of $X$  (cf. \cite[Theorem 1.4.10]{brion-kumar}). By abuse of language, we say that $X$ is F-split by $\si\in H^0(X,\omega_X^{-1})$. 

The latter category includes \emph{spherical varieties} (cf. \cite{br-in}), in particular projective homogeneous varieties and toric varieties, and also \emph{Fano varieties} (cf. \cite[Exercise 1.6.E(5)]{brion-kumar}). 
\end{m-remark}


\end{document}